\def\draft{n}
\documentclass{amsart}
\usepackage[headings]{fullpage}
\usepackage{amssymb,epic,eepic,epsfig,amsbsy,amsmath}
\usepackage[mathscr]{euscript}

\theoremstyle{plain}

\newtheorem{theorem}{Theorem}
\newtheorem{proposition}{Proposition}[section]
\newtheorem{lemma}[proposition]{Lemma}

\theoremstyle{definition}
\newtheorem{definition}[proposition]{Definition}

\newtheorem{question}{Question}

\theoremstyle{remark}

\newtheorem{remark}[proposition]{Remark}

\def\printname#1{
    \if\draft y
        \smash{\makebox[0pt]{\hspace{-0.5in}
            \raisebox{8pt}{\tt\tiny #1}}}
    \fi
}

\newlength{\standardunitlength}
\catcode`\@=11
\long\def\@makecaption#1#2{%
     \vskip 10pt

\setbox\@tempboxa\hbox{
       \small\sf{\bfcaptionfont #1. }\ignorespaces #2}%
     \ifdim \wd\@tempboxa >\captionwidth {%
         \rightskip=\@captionmargin\leftskip=\@captionmargin
         \unhbox\@tempboxa\par}%
       \else
         \hbox to\hsize{\hfil\box\@tempboxa\hfil}%
     \fi}
\font\bfcaptionfont=cmssbx10 scaled \magstephalf
\newdimen\@captionmargin\@captionmargin=2\parindent
\newdimen\captionwidth\captionwidth=\hsize
\catcode`\@=12

\def\lbl#1{\label{#1}\printname{#1}}

\def\BN{\mathbb N}
\def\BZ{\mathbb Z}

\def\BQ{\mathbb Q}

\def\BC{\mathbb C}

\def\D{\Delta}

\def\M{\mathcal M}

\def\calR{\mathcal R}

\def\a{\alpha}

\def\l{\lambda}

\def\S{\Sigma}

\def\la{\langle}
\def\ra{\rangle}

\def\b{\beta}

\newcommand{\Ker}{\operatorname{Ker}}

\def\longto{\longrightarrow}

\def\pt{\partial}

\def\Tr{\mathrm{Tr}}

\def\calL{\mathcal{L}}

\def\alg{\mathrm{alg}}
\def\rat{\mathrm{rat}}
\def\hol{\mathrm{hol}}
\def\conv{\mathrm{conv}}

\begin{document}

\title[Algebraic $G$-functions associated to matrices over a group-ring]{
Algebraic $G$-functions associated to matrices over a group-ring}

\author{Stavros Garoufalidis}
\address{School of Mathematics \\
         Georgia Institute of Technology \\
         Atlanta, GA 30332-0160, USA \\
         {\tt http://www.math.gatech} \newline {\tt .edu/$\sim$stavros } }
\email{stavros@math.gatech.edu}
\author{Jean Bellissard}
\address{School of Mathematics \\
         Georgia Institute of Technology \\
         Atlanta, GA 30332-0160, USA \\
         {\tt http://www.math.gatech} \newline {\tt .edu/$\sim$jeanbel } }
\email{jeanbel@math.gatech.edu}

\thanks{The authors were supported in part by NSF. \\
\newline
1991 {\em Mathematics Classification.} Primary 57N10. Secondary 57M25.
}

\date{October 1, 2007}


\begin{abstract}
Given a square matrix with elements in the group-ring of a group, one
can consider the sequence formed by the trace (in the sense of the group-ring)
of its powers. We prove that the corresponding generating series is an
algebraic $G$-function (in the sense of Siegel) when the group is free of 
finite rank. Consequently, it follows that the norm of such elements is
an exactly computable algebraic number, and their Green function is
algebraic. Our proof uses the notion of rational and algebraic power series 
in non-commuting variables and is an easy application of a theorem of Haiman. 
Haiman's theorem uses results of linguistics regarding regular and 
context-free language. On the other hand, when the group is free abelian of 
finite rank, then the corresponding generating series is a $G$-function. 
We ask whether the latter holds for general hyperbolic groups. 
\end{abstract}

\maketitle

\tableofcontents

\section{Introduction}
\lbl{sec.intro}

\subsection{Algebricity of the Green's function for the free group}
\lbl{sub.algfree}

Given a group $G$, consider the group-algebra 
$\BQ[G]$, and define a {\em trace} map:
\begin{equation}
\lbl{eq.trace}
\Tr: \BQ[G] \longto \BC, \qquad
\Tr(P)=\text{constant term of}\,\, P
\end{equation}
where the constant term is the coefficient of the identity element of $G$.
Let $M_N(R)$ denote the set of $N$ by $N$ matrices with entries in a ring $R$.
We can extend the trace to the algebra $M_N(\BQ[G])$ by:
\begin{equation}
\lbl{eq.traceM}
\Tr: M_N(\BQ[G]) \longto \BC, \qquad
\Tr(P)= \sum_{j=1}^N \Tr(P_{jj}).
\end{equation}

\begin{definition}
\lbl{def.moments}
Given $P \in M_N(\BQ[G])$, 
consider the sequence $(a_{P,n})$ 
\begin{equation}
\lbl{eq.aMn}
a_{P,n}=\Tr(P^n)
\end{equation}
and the generating series
\begin{equation}
\lbl{eq.Rz}
R_P(z)=\sum_{n=0}^\infty a_{P,n} z^n. 
\end{equation}
\end{definition}

Let $F_r$ denote the free group of rank $r$.

\begin{theorem}
\lbl{thm.1}
The Green's function $R_P(z)$ of every element $P$ of $M_N(\BQ[F_r])$ is 
algebraic. 
\end{theorem}

Theorem \ref{thm.1} appears in the cross-roads of several areas of research:

\begin{itemize}
\item[(a)] 
operator algebras
\item[(b)]
free probability
\item[(c)]
linguistics and context-free languages
\item[(d)]
non-commutative combinatorics
\item[(e)]
mathematical physics
\end{itemize}

In fact, Woess proves Theorem \ref{thm.1} when $N=1$ using linguistics
and context-free languages; see \cite{Wo1,Wo2}. Voiculescu proves Theorem 
\ref{thm.1} using the $R$ and $S$ transforms of free probability; see
\cite{Vo1,Vo2}. For additional results using free probability, see
\cite{Ao,CV} and also \cite{Le1,Le2}.

It is well-known that Theorem \ref{thm.1} provides an exact calculation
of the norm of $P \in M_N(\BQ[F_r]) \subset  M_N(L(F_r))$, where
$L(F_r)$ denotes the {\em reduced $\BC^*$-algebra} completion of the 
group-algebra  $\BC[F_r]$. For a detailed discussion, see the above references.

Our proof of Theorem \ref{thm.1} uses the notion of an algebraic function in 
non-commuting variables and a theorem of Haiman, which itself is based on
a theorem of Chomsky-Sch\"utzenberger on context-free languages. A by-product
of our proof is the fact that the moment generating series is a matrix
of algebraic power series in non-commuting variables (see Proposition
\ref{prop.alg}), which is a statement a priori
stronger than Theorem \ref{thm.1}.

An alternative proof of Theorem \ref{thm.1} uses methods from functional 
analysis, and most notably the {\em Schur complement method} (see below). 
We will discuss in detail the first proof and postpone the third proof to a 
later publication. Either proof explains the close relation between the 
differential properties
of the generating function $R_P(z)$ and the word problem in $G$.

\section{The case of the free abelian group}
\lbl{sec.freeab}

\subsection{Holonomic, algebraic and $G$-functions}
\lbl{sub.holonomic}

A priori, $R_P(z)$ is only a formal power series. However, it is easy to
see that $(a_{P,n})$ is bounded exponentially by $n$, which implies that
$R_P(z)$ defines an analytic function in a neighborhood of $z=0$.
The paper is concerned with differential/algebraic properties of 
the function $R_P(z)$. 
Algebraic and holonomic functions are well-studied objects. Let us
recall their definition here.

\begin{definition}
\lbl{def.holo}
\rm{(a)} A {\em holonomic} function $f(z)$ is one that
satisfies a linear differential equation with polynomial 
coefficients. In other words, we have:
$$
c_d(z) f^{(d)}(z) + \dots + c_0(z) f(z)=0
$$
where $c_j(z) \in \BQ[z]$ for all $j=0,\dots, d$
and $f^{(j)}(z)=d^j/dz^j f(z)$. 
\newline
\rm{(b)} An {\em algebraic} function $f(z)$ is one
that satisfies a polynomial equation:
$$
Q(f(z),z)=0
$$
where $Q(y,z) \in \BQ[y,z]$.
\end{definition}

Lesser known to the combinatorics community are $G$-functions, which 
originated in the work Siegel on arithmetic problems in 
elliptic integrals, and transcendence problems in number theory;
see \cite{Si}. Holonomic $G$-functions originate naturally in 
\begin{itemize}
\item[(a)]
algebraic geometry, related to the regularity properties of the Gauss-Manin 
connection, see for example \cite{De,Ka,Ma},
\item[(b)]
arithmetic, see for example \cite{An2,Bm,DGS},
\item[(c)]
enumerative combinatorics, as was recently shown in \cite{Ga2}.
\end{itemize}

\begin{definition}
\lbl{def.Gfunctions}
A {\em $G$-function} $f(z)=\sum_{n=0}^\infty a_n z^n$ is one 
which satisfies the following conditions:
\begin{itemize}
\item[(a)]
for every $n \in \BN$, we have $a_n \in \overline{\BQ}$, 
\item[(b)]
there exist a constant $C_f>0$ such that for every $n \in \BN$
we have: $|a_n| \leq C_f^n$ (for every conjugate of $a_n$)
and the common denominator of $a_0,\dots,a_n$
is less than or equal to $C_f^n$. 
\item[(c)] $f(z)$ is holonomic.
\end{itemize}
\end{definition}

The next theorem summarizes the analytic continuation and the shape of 
the singularities of algebraic functions and $G$-functions. Part (a) follows
from the general theory of differential equations (see eg. \cite{Wa}), 
parts (b) and (d) follow from \cite[Lem.2.2]{CSTU} 
(see also \cite{DGS} and \cite{DvdP})
and (c) follows from a combination of Katz's theorem, Chudnovsky's theorem
and Andr\'e's theorem; see \cite[p.706]{An2} and also \cite{C-L}.

\begin{theorem}
\lbl{thm.0}
\rm{(a)}
A holonomic function $f(z)$ can be analytically continued as a multivalued
function in $\BC\setminus\S_f$ where $\S_f \subset \overline{\BQ}$ 
is the finite set of singular points of $f(z)$. 
\newline
\rm{(b)}
Every algebraic function $f(z)$ is a $G$-function.
\newline
\rm{(c)}
In a neighborhood of a singular point $\l \in \S_f$, a $G$-function 
$f(z)$ can be written as a finite sum of germs of the form:
\begin{equation}
\lbl{eq.nilsson}
(z-\l)^{\a_{\l}}(\log(z-\l))^{\b_{\l}} h_{\l}(z-\l)
\end{equation}
where $\a_{\l} \in \BQ$, $\b_{\l} \in \BN$, and $h_{\l}$ a holonomic 
$G$-function.
\newline
\rm{(d)}
In addition, $\b_{\l}=0$ if $f(z)$ is algebraic. 
\end{theorem}

\begin{remark}
\lbl{rem.Gfunctions}
Local expansions of the form \eqref{eq.nilsson} are known in the literature as
Nilsson series (see \cite{Ni}), and minimal order linear differential 
equations that they satisfy are known to be regular singular, with rational
exponents $\{a_{\l}\}$ and quasi-unipotent monodromy. For a discussion,
see \cite{Ka,Ma,Ga2} and references therein.
\end{remark}

It is classical and easy to show that the existence of analytic continuation
of a function implies the existence of asymptotic expansion of its Taylor
series; see for example \cite{Ju, Co} and also \cite[Sec.7]{CG} and \cite{Ga2}.

\begin{lemma}
\lbl{lem.anas}
If $f(z)=\sum_{n=0}^\infty a_n z^n$ is holonomic and analytic at $z=0$, then
the $n$th Taylor coefficient $a_n$ has an asymptotic expansion
in the sense of Poincar\'e
\begin{equation}
\lbl{eq.trans}
a_n \sim \sum_{\l \in \S} \l^{-n} n^{-\a_{\l}-1} (\log n)^{\b_{\l}} 
\sum_{s=0}^\infty \frac{c_{\l,s}}{n^s}
\end{equation}
where $\S_f$ is the set of singularities of $f$, $\a_{\l}, \b_{\l} \in \BQ$,
and $c_{\l,s} \in \BC$.
\end{lemma}

\subsection{The case of the free abelian group}
\lbl{sub.invcom}

In this section we will summarize what is known about the generating
functions $R_P(z)$ when $G=\BZ^r$ is the free abelian group or rank $r$.
The next theorem is shown in \cite{Ga2}, using Andre\'e main theorems from
\cite{An2}. An alternative proof uses the regular holonomicity of the 
Gauss-Manin connection and the rationality of its exponents.
This was kindly communicated to us by C. Sabbah (see also \cite{DvK}). 
Holonomicity of $R_P(z)$ also follows from a fundamental result of
Wilf-Zeilberger, explained in \cite{Ga2}.

\begin{theorem}
\lbl{thm.2}\cite{Ga2}
For every $P \in M_N(\BQ[\BZ^r])$, $R_P(z)$ is a $G$-function.
\end{theorem}

\subsection{A complexity remark}
\lbl{sub.complexity}

Given $P \in M_N(\BQ[F_r])$ (resp. $P \in M_N(\BQ[F_r])$), one may ask for the 
complexity of a minimal polynomial $Q(y,z) \in \BQ[y,z]$ 
(resp. minimal degree differential operator 
$D(z,\pt_z) \in \BQ \la z, \pt z \ra$) so that $Q(R_P(z),z)=0$ 
(resp. $D(z,\pt_z) R_P(z)=0$). One expects that the $y$-degree of $Q(y,z)$
and the $\pt_z$-degree of $D(z,\pt_z)$ is exponential in the {\em complexity}
of $P$, where the latter can be defined to be the degree of $P$ and 
the maximum of the absolute values of the coefficients of the enrties of $P$.
This prohibits explicit calculations in general.

\subsection{Acknowledgement}

The second author wishes to thank R. Gilman, F. Flajolet, L. Mosher, C. Sabbah
and D. Zeilberger for stimulating conversations and D. Voiculescu for 
bringing into the attention relevant literature on free probability.

\section{A theorem of Haiman and a proof of Theorem \ref{thm.1}}
\lbl{sec.haiman}

In \cite{Ha} Haiman proves the following theorem.

\begin{theorem}
\lbl{thm.ha}\cite{Ha}
Let $K$ be a field with a rank 1 discrete valuation $v$; $K_v$ its
completion with respect to the metric induced by $v$. Let $f(x_1,\dots,x_r,
y_1,\dots,y_r)$ be a rational power series over $K$ in non-commuting
indeterminants. Any coefficient of $f(x_1,\dots,x_r,x_1^{-1},\dots, x_r^{-1})$
converging over $K_v$ is algebraic over $K$.
\end{theorem}

Letting $K=\BQ(z)$, and $K_v=\BQ((z))$ the ring of formal Laurent series in
$z$, and considering the element $(1-zP)^{-1}$, where $P \in \M_N(\BQ[F_r])$,
gives an immediate proof of Theorem \ref{thm.1}.

In the next section we will give a detailed description of Haiman's
argument which exhibits a close relation to linguistics, as well as
an obstruction to generalizing Theorem \ref{thm.1} to groups other than
the free group.

\section{Algebraic and rational functions in noncommuting variables}
\lbl{sec.noncommuting}

\subsection{Rational, algebraic and holonomic functions in one variable}
\lbl{sub.1var}

In
this section all functions will be analytic in a neighborhood of $z=0$.
Let $\BQ^{\rat}_0(z)$, $\BQ^{\alg}_0(z)$ and $\BQ^{\hol}_0(z)$
denote respectively the set of rational, algebraic and holonomic
functions, analytic at $z=0$. Let $\BQ[[z]]$ denote the set of formal
power series in $z$. Using the injective Taylor series map around $z=0$, 
we will consider $\BQ^{\rat}_0(z)$, $\BQ^{\alg}_0(z)$ and $\BQ^{\hol}_0(z)$
as subsets of $\BQ[[z]]$:

\begin{equation}
\lbl{eq.inclusion}
\BQ^{\rat}_0(z) \subset \BQ^{\alg}_0(z) \subset \BQ^{\hol}_0(z) 
\subset \BQ[[z]].
\end{equation}
$\BQ[[z]]$ has {\em two} multiplications:
\begin{itemize}
\item
the usual multiplication of formal power series
\begin{equation}
\lbl{eq.usualm}
(\sum_{n=0}^\infty a_n z^n ) \cdot (\sum_{n=0}^\infty b_n z^n )=
\sum_{n=0}^\infty (\sum_{k=0}^n a_k b_{n-k}) z^n .
\end{equation}
\item
The {\em Hadamard product}:
\begin{equation}
\lbl{eq.hadamardm}
(\sum_{n=0}^\infty a_n z^n ) \circledast (\sum_{n=0}^\infty b_n z^n )=
\sum_{n=0}^\infty a_n b_n z^n .
\end{equation}
\end{itemize}
With respect to the usual multiplication, $\BQ[[z]]$ 
is an algebra and $\BQ^{\rat}_0(z)$, $\BQ^{\alg}_0(z)$ and $\BQ^{\hol}_0(z)$
are subalgebras.
In case two power series are convergent in a neighborhood of zero, so
is their Hadamard product. Hadamard, Borel and Jungen studied the 
analytic continuation and the singularities of the Hadamard product of two 
functions; see \cite{Bo, Ju}. Their method used an integral representation
of the Hadamard product, and a deformation of the contour of integration;
see \cite[Fig.2,p.303]{Ju}.
Let us summarize these classical results.

\begin{theorem}
\lbl{thm.1var}
(a) If $f$ and $g$ are rational, so is $f \circledast g$.
\newline
(b) If $f$ is rational and $g$ is algebraic, then $f \circledast g$ is
algebraic.
\newline
(c) If $f$ and $g$ are holonomic (resp. regular holonomic with rational
exponents),  so is $f \circledast g$.
\newline
(d) If $f$ and $g$ are algebraic, then $f \circledast g$ is not necessarily
algebraic.
\end{theorem}

For a proof, see Thm.7, Thm.8, Theorem E and
the example of p.298 from \cite{Ju}.

\subsection{Rational and algebraic functions in noncommuting variables}
\lbl{sub.noncommute}

In this section we discuss a generalization of the previous section
to non-commuting variables. Let $X$ be a finite set, and let $X^*$ denote
the free monoid on $X$. In other words, $X$ consists of 
the set of all words in $X$, including the empty word $e$.
Let $\BQ\la X \ra$ (resp. $\BQ\la\la X \ra\ra$) denote the algebra of 
polynomials (resp. formal power series) in non-commuting variables.
In \cite{Sh}, Sch\"utzenberger defines the notion of a {\em rational} and
an {\em algebraic} power series in non-commuting variables. Let 
$\BQ^{\rat}\la X\ra$ and $\BQ^{\alg}\la X\ra$ denote the sets of rational
(resp. algebraic) power series. Then, we have an inclusion:

\begin{equation}
\lbl{eq.inclusion2}
\BQ^{\rat}\la X\ra \subset \BQ^{\alg}\la X\ra \subset \BQ\la\la X \ra\ra.
\end{equation}
$\BQ\la\la X \ra\ra$ has two multiplications: 

\begin{itemize}
\item
the usual multiplication of formal power series in non-commuting variables:
\begin{equation}
\lbl{eq.usualmX}
(\sum_{w \in X^*} a_w w ) \cdot (\sum_{w \in X^*} b_w w)=
\sum_{w \in X^*} (\sum_{w',w'': w'w''=w} a_{w'} b_{w''}) w .
\end{equation}
\item
The Hadamard product:
\begin{equation}
\lbl{eq.hadamardmX}
(\sum_{w \in X^*} a_w w ) \circledast (\sum_{w \in X^*} b_w w)=
\sum_{w \in X^*}  a_w b_w w .
\end{equation}
\end{itemize}
With respect to the usual multiplication, $\BQ\la\la X \ra\ra$ is a 
non-commutative algebra and $\BQ^{\rat}\la X\ra$ and $\BQ^{\alg}\la X\ra$ are
subalgebras. We have  the following analogue of Theorem 
\ref{thm.1var}.

\begin{theorem}
\lbl{thm.Xvar}
\cite[Pro.2.2]{Sh} 
(a) If $f \in \BQ^{\rat}\la X\ra $ and $g \in \BQ^{\rat}\la X\ra$, then
$f \circledast g \in \BQ^{\rat}\la X\ra$.
\newline
(b) If $f \in \BQ^{\rat}\la X\ra$ and $g \in \BQ^{\alg}\la X\ra$, then 
$f \circledast g \in \BQ^{\alg}\la X\ra$.
\end{theorem}

\begin{remark}
\lbl{rem.changering}
The notion of rational and algebraic functions works for an arbitrary ring
$\calR$ of characteristic zero, instead of $\BQ$. Theorem \ref{thm.Xvar}
is still valid.
\end{remark}

\subsection{Proof of Theorem \ref{thm.1}}
\lbl{sub.thm1}

Let $F_r$ denote the free group of rank $r$ with generating set 
$\{u_1,\dots,u_r\}$, and 
$$
X=\{x_1,\dots,x_r,\overline{x}_1,\dots,\overline{x}_r\}.
$$ 
Consider the monoid map:

\begin{equation}
\lbl{eq.pi}
\pi: X^* \longto F_r, \qquad \pi(x_i)=u_i, \qquad \pi(\overline{x}_i)=u_i^{-1}.
\end{equation}
The kernel $\Ker(\pi)$ of $\pi$ is the set of those words in $X$  
which reduce to the identity under the relations 
$x_i\overline{x}_i=\overline{x}_ix_i=e$.
Let 
\begin{equation}
\lbl{eq.delta}
\D=\sum_{w \in \Ker(\pi)} w \,\, \in \BQ\la\la X \ra\ra.
\end{equation}

\noindent
The next proposition is attributed to Chomsky-Sch\"utzenberger by Haiman.
For a proof, see \cite[Sec.3]{Ha}.

\begin{proposition}
\lbl{prop.CS}
\cite{CS}
$\D$ is algebraic.
\end{proposition}

The map $\pi$ has a right inverse (that satisfies $\pi \circ \iota=I_{F_r}$):

\begin{equation}
\lbl{eq.invpi}
\iota: F_r \longto X
\end{equation}
defined by mapping a reduced word in $u_i$ to a corresponding word in $X$.
For every $f \in \BQ[F_r]$ we have a key relation between 
trace and Hadamard product:

\begin{equation}
\lbl{eq.key}
\Tr(f)=\phi(\iota(f) \circledast \D)
\end{equation}
where $\phi$ is a $\BQ$-linear map defined by:

\begin{equation}
\lbl{eq.phi}
\phi: \BQ\la X\ra \longto \BQ, \qquad \phi(w)=1 
\qquad \text{for} \,\, w \in X^*.
\end{equation}

Now, fix $P \in M_N(\BQ[F_r])$. 
Let $\D_N$ denote the $N$ by $N$ matrix with entries equal to $\D$,
and $\calR=\BQ(z)$. Let
$$
P_z=z\iota(P) \in M_N(\calR\la X \ra), \qquad
P_z^*=\sum_{n=0}^\infty P_z^n \in M_N(\calR\la\la X \ra\ra).
$$
Notice that $P_z^*$ is well-defined since $P_z$ has
no $z$-constant term.

\begin{lemma}
\lbl{lem.rat}
We have:
\begin{equation}
\lbl{eq.Pz}
P_z^* \in M_N(\calR^{\rat}\la X \ra).
\end{equation}
\end{lemma}

\begin{proof}
$P_z^*$ satisfies the matrix equation
$$
(1-P_z) P_z^*=I
$$
with entries in $\calR\la X \ra$.
\end{proof}

Lemma \ref{lem.rat}, together with Propositions \ref{prop.CS} and part (b) of
\ref{thm.Xvar} imply the following result, which we can think as a 
noncommutative analogue of Theorem \ref{thm.1}.

\begin{proposition}
\lbl{prop.alg}
For every $P \in M_N(\BQ[F_r])$, we have:
\begin{equation}
\lbl{eq.alg1}
\sum_{n=0}^\infty z^n (\iota(P))^n \circledast \D_N \in 
M_N(\calR^{\alg}\la X \ra).
\end{equation}
\end{proposition}

Consider the abelianization ring homomorphism:
\begin{equation}
\lbl{eq.psi}
\psi: \calR\la\la X\ra\ra \longto \calR[[X]]
\end{equation}
where $\calR[[X]]$ is the formal power series ring in commuting variables.
Haiman proves the following:

\begin{proposition}
\lbl{prop.ha2}
\cite[Prop.3.3]{Ha}
If $f \in \calR^{\alg}\la X\ra$, then $\psi(f)$ is algebraic over
$\calR(X)$.
\end{proposition}

It follows that $\psi(P_z^* \circledast \D_N ) \in M_N(\calR^{\alg}(X))$. 
Consider now the subalgebra $\calR^{\conv}[[X]]$ of $\calR[[X]]$
that contains all elements of the form

$$
\sum_{w \in X^*} a_w w 
$$
where $a_w \in z^{l(w)}\BQ[[z]]$, where $l(w)$ denotes the length of $w$.
Then, we can define an algebra map:

\begin{equation}
\lbl{eq.phiz}
\phi_z: \calR^{\conv}[[X]] \longto \BQ[[z]], \qquad \phi_z(w)=1 \,\,
\qquad \text{for} \,\, x \in X.
\end{equation}

Haiman shows that if $f \in \calR^{\alg}(X) \cap \calR^{\conv}[[X]]$, then
$\phi_z(f) \in \BQ^{\alg}$.
To state our final conclusion, we define for $1 \leq i,j \leq N$, the
sequence $(a^{ij}_{P,n})$ by 
$$
a^{ij}_{P,n}=\Tr((P^n)_{ij})
$$
and the matrix of generating series $A_P(z) \in M_N(\BQ[[z]])$ by:
$$
(A_P(z))_{ij}=\sum_{n=0}^\infty a^{ij}_{P,n} z^n.
$$

\begin{lemma}
\lbl{lem.compute}
We have:
\begin{equation}
\lbl{eq.PwD}
(\phi_z \circ \psi)( P_z^* \circledast \D_N)= A_P(z).
\end{equation}
Thus, $A_P(z) \in M_N(\BQ^{\alg}_0(z))$.
\end{lemma}

\begin{proof}
Equation \eqref{eq.PwD} follows from Equation \eqref{eq.key}.
The conclusion follows from the above discussion.
\end{proof}

\noindent
Thus, the entries of $A_P(z)$ are algebraic functions, convergent at $z=0$. 
Since by definition we have:
$$
R_P(z)=\sum_{i=1}^N (A_P(z))_{ii}
$$
it follows that $R_P(z) \in \BQ^{\alg}_0(z)$. This completes the
proof of Theorem \ref{thm.1}.
\qed

\section{Some Linguistics}
\lbl{sec.linguistics}

\subsection{Regular and context-free languages}
\lbl{sub.ling}

Haiman's proof uses the key Proposition \ref{prop.CS} from linguistics.
Let us recall some concepts from this field. See for example \cite{BR,Li,Ya}
and references therein. Given a finite set $X$ (the
alphabet), a language $L$ is a collection of words in $X$. In other words,
$\calL \subset X^*$. The {\em generating series} $F_L$ of a language is:
$$
F_{\calL}=\sum_{w \in L} w \in \BQ\la\la X \ra\ra.
$$
It follows that for two languages $\calL_1$ and $\calL_2$ we have:
$$
F_{\calL_1 \cap \calL_2}=F_{\calL_1} \circledast F_{\calL_2} .
$$
A language $L$ is called {\em rational} (resp. {\em context-free})
iff $F_L \in \BQ^{\rat}(X)$ (resp. $F_L \in \BQ^{\alg}(X)$).
In this context, Theorem \ref{thm.Xvar} takes the following form:

\begin{theorem}
\lbl{thm.ling}\cite{CS}
(a) If $\calL_1$ and $\calL_2$ are rational languages, so is 
$\calL_1 \cap \calL_2$.
\newline
(b) If $\calL_1$ is rational and $\calL_2$ is (unambiguous)
context-free, then $\calL_1 \cap \calL_2$ is (unambiguous) context-free.
\end{theorem}

It was pointed out to us independently by D. Zeilberger and F. Flajolet that 
the above theorem essentially proves Theorem \ref{thm.1}.

\subsection{Some questions}
\lbl{sub.question}

Let us end this short paper with some questions.
Despite the similarity in their statements and the multitude of proofs, 
Theorems \ref{thm.1} and \ref{thm.2} have 
different assumptions, different proofs and different conclusions. 

Consider a generating set $X$ for a group $G$ such that every element of $G$
can be written as a word in $X$ with nonnegative exponents.
Given $X$ and $G$, let $\calL_X$ denote the set of all words in $X$ that map
to the identity in $G$. Deciding membership in $\calL_X$ is the {\em 
word problem} in $G$.

\begin{definition}
\lbl{def.wordproblem}
A group $G$ has {\em context-free word problem} if it has a
generating set $X$ such that the language $\calL_X$ is {\em context-free}.
\end{definition}

The proof of Theorem \ref{thm.1} applies to groups with
a context-free word problem. Miller-Schupp classified those groups.
In \cite{MS} Miller-Schupp prove that $G$ has context-free word problem
iff $G$ has a free finite-index subgroup.

On the other hand, if $G$ is the fundamental group of a hyperbolic manifold
of dimension not equal to $2$, then $G$ does not have a free finite-index
subgroup.

Thus, the linguistics proof of Theorem \ref{thm.1} does not apply to the
case of hyperbolic groups in dimension three. Neither does it apply
to the case of $\BZ^r$ since the latter does not have context-free word 
problem. 

\begin{question}
\lbl{que.1}
If $P$ is a hyperbolic group and $P \in M_N(\BQ[G])$, is it true that
$R_P(z)$ is a $G$-function?
\end{question}

The question may be relevant to low dimensional topology, 
when one tries to compute the 
$\ell^2$-torsion of a hyperbolic manifold using Luecke's theorem; \cite{Lu}.
In that case, the matrix $P$ comes from Fox (free differential) calculus
of a presentation of the fundamental group $G$ of the hyperbolic manifold.
See also \cite{DL}.

\begin{question}
\lbl{que.2}
Given $P \in M_N(\BQ[F_r])$, consider the abelianization $P^{\mathrm{ab}} 
\in M_N(\BQ[\BZ^r])$, and the $G$-functions $R_P(z)$ and 
$R_{P^{\mathrm{ab}}}(z)$. How are the singularities of $R_P(z)$ and 
$R_{P^{\mathrm{ab}}}(z)$ related?
\end{question}

\begin{question}
\lbl{que.3}
What is a holonomic function in non-commuting variables?
\end{question}

\section{A functional analysis interpretation of Theorem \ref{thm.1}}
\lbl{sub.physics}

The present paper is focusing on results and techniques inspired by 
algebra, non-commutative algebraic combinatorics. However it is worth 
mentioning that Theorem \ref{thm.1} has applications to problems coming from 
functional analysis, spectral theory, and the spectrum of
Schr\"odinger operators.  For instance, the Schr\"odinger equation 
describing the electron motion in a $d$-dimensional periodic crystal, can be 
well approximated by the difference equation on a lattice of same dimension. 
The corresponding operator can be seen as an element of the group ring of 
$\BZ^d$. The function $R_P(z)$ defined previously is noting but the diagonal 
element of the resolvent and is used to compute the spectral measure, 
through the Charles de la Vall\'ee Poussin theorem. 
There are instances for which, this operator is better approximated by the 
free group analog. For instance the {\em retracable path approximation} was 
used by Brinkman and Rice \cite{BR71} in 1971 to treat the effect of 
spin-orbit coupling in the Hall effect, while it was used in \cite{BFZ} to 
compute the electronic Density of States when the electron is submitted to 
a random magnetic field.
The same operator, seen as an element of the free group ring, is used to 
describe various infinite dimension approximations. The seminal work of 
Georges and Kotliar \cite{GK} used this free group approximation to give 
the first model known with a {\em Mott-Hubbard} transition.

Another domain in which the Theorem \ref{thm.1} may apply is the 
Voiculescu Theory of {\em Free Probability} \cite{Vo2,Vo3}. The so-called 
{\em $R$-transform} used to treat the convolution of free random variables, 
is also based upon the Schur complement formula. In particular the free 
central limit theorem asserts that a sum of identically distributed free 
random variable obey the semicircle law, is a special case of the present 
result. 

Besides the two proofs of Theorem \ref{thm.1} discussed 
in this paper, the algebraic 
character of $R_P(z)$ can also be deduced from the used of the Schur 
complement method \cite{Schur17}. This is what makes the free group 
approximation so attractive to theoretical physicists. 
This method, also known under the name of {\em Feshbach method}
\cite{Feshbach58} is used in many domains of Physics, Quantum Chemistry, 
Solid State Physics, Nuclear Physics, to reduce the Hilbert space to a 
finite dimensional one and make the problem amenable to numerical 
calculations. However, very few Mathematical Physicists have paid attention 
to the fact that algebraicity or holonomy can give rise to results 
concerning the explicit computation of the spectral radius, or more 
generally, to the band edges, of the Hamiltonian they consider. This 
later problem is known to be notably hard with other methods.

For the benefit of the reader, we include some history of that
method. The Schur complement method \cite{Schur17} is widely used in 
numerical analysis under this name, while Mathematical Physicists prefer 
the reference to Feshbach \cite{Feshbach58}. In Quantum Chemistry, the common
reference is Feshbach-Fano \cite{Fano35} or Feshbach-L\"owdin 
\cite{Lowdin62}. This method is used in various
algorithms in Quantum Chemistry ({\em ab initio} calculations), in Solid 
State Physics
(the muffin tin approximation, LMTO) as well as in Nuclear Physics. The 
formula used
above is found in the original paper of Schur \cite[p.217]{Schur17}. 

The formula has been proposed also by an astronomer Tadeusz Banachiewicz 
in 1937, even though closely related results were obtained in 1923 by Hans 
Boltz and in 1933 by Ralf Rohan \cite{PS04}. Applied to the Green function 
of a selfadjoint operator with finite rank perturbation, it becomes the 
Kre\u{\i}n formula \cite{Krein46}.

Let us end this section with a small dictionary that
compares our notions with those in physics.

\begin{center}
\begin{tabular}{|l|l|}
\hline
$H \in M_N(\BQ[F_r])$ & Hamiltonian \\ \hline
$1/(z-H)$ & resolvant  \\ \hline
$1/z R_H(1/z)$ & trace of the resolvant \\ \hline
$\Tr(H^n)$ & $n$th moment of $H$ \\ \hline
\end{tabular}
\end{center}

\ifx\undefined\bysame
    \newcommand{\bysame}{\leavevmode\hbox
to3em{\hrulefill}\,}
\fi

\end{document}